\documentclass[10pt,letterpaper,twoside,reqno]{amsart}

\makeatletter{}
\usepackage{fixltx2e}                       
\usepackage[usenames,dvipsnames]{xcolor}
\usepackage{fancyhdr}
\usepackage{amsmath,amsfonts,amsbsy,amsgen,amscd,mathrsfs,amssymb,amsthm}
\usepackage{subfig}
\usepackage{url}

\usepackage{mathtools}

\mathtoolsset{showonlyrefs}

\usepackage[font=small,margin=0.25in,labelfont={sc},labelsep={colon}]{caption}

\usepackage{tikz}
\usepackage{microtype}
\usepackage{enumitem}

\definecolor{dark-gray}{gray}{0.3}
\definecolor{dkgray}{rgb}{.4,.4,.4}
\definecolor{dkblue}{rgb}{0,0,.5}
\definecolor{medblue}{rgb}{0,0,.75}
\definecolor{rust}{rgb}{0.5,0.1,0.1}

\usepackage[colorlinks=true]{hyperref}

\hypersetup{urlcolor=rust}
\hypersetup{citecolor=dkblue}
\hypersetup{linkcolor=dkblue}

\usepackage{setspace}

\usepackage{graphicx}
\usepackage{booktabs,longtable,tabu} 
\usepackage{multirow} 
\usepackage{float}

\usepackage{fourier}

\usepackage{bm} 
\graphicspath{{figures/}}

\newtheorem{theorem}{Theorem}[section]

\theoremstyle{definition}

\newcommand{\term}{\emph}

\numberwithin{equation}{section} 
\numberwithin{figure}{section}
\numberwithin{table}{section}

\floatstyle{plaintop}
\newfloat{recipe}{thp}{lor}
\floatname{recipe}{Recipe}
\numberwithin{recipe}{section}

\providecommand{\mathbold}[1]{\bm{#1}}

\renewcommand{\phi}{\varphi}

\newcommand{\half}{\tfrac{1}{2}}

\newcommand{\Id}{\mathbf{I}}

\providecommand{\mathbbm}{\mathbb} 
\newcommand{\R}{\mathbbm{R}}

\newcommand{\N}{\mathbbm{N}}
\newcommand{\Z}{\mathbbm{Z}}
\newcommand{\Sym}{\mathbb{H}}

\newcommand{\abs}[1]{\left\vert {#1} \right\vert}

\newcommand{\Expect}{\operatorname{\mathbb{E}}}

\DeclareMathOperator{\Var}{Var}

\newcommand{\vct}[1]{\mathbold{#1}}
\newcommand{\mtx}[1]{\mathbold{#1}}

\newcommand{\adj}{*}

\newcommand{\trace}{\operatorname{tr}}

\newcommand{\cone}{\operatorname{cone}}

\title{Integer Factorization of a Positive-Definite Matrix}

\author[J.~A.~Tropp]{Joel~A.~Tropp}

\date{26 May 2015.  Revised: 31 May 2015, 1 July 2015, and 3 August 2015.}


\begin{document}

\begin{abstract}
This paper establishes that every positive-definite matrix
can be written as a positive linear combination of outer products of
integer-valued vectors whose entries are bounded by the geometric
mean of the condition number and the dimension of the matrix.
\end{abstract}

\maketitle

\section{Motivation}

This paper addresses a geometric question that arises
in the theory of discrete normal approximation~\cite{BLX15:Multivariate-Approximation}
and in the analysis of hardware for implementing matrix multiplication~\cite{LUW15:Factorization-Analog}.
The problem requires us to represent a nonsingular
covariance matrix as a positive linear combination of
outer products of integer vectors.  The theoretical challenge is to
obtain an optimal bound on the magnitude of the integers
required as a function of the condition number of the
matrix.  We establish the following result.

\begin{theorem} \label{thm:main}
For positive integers $m$ and $d$, define a set of bounded integer vectors:
$$
\Z^d_m := \big\{ \vct{z} \in \Z^d : \abs{z_i} \leq m \text{ for $i=1,\dots,d$} \big\}.
$$
Let $\mtx{A}$ be a real $d \times d$ positive-definite matrix with
(finite) spectral condition number
$$
\kappa(\mtx{A}) := \lambda_{\max}(\mtx{A}) / \lambda_{\min}(\mtx{A}),
$$
where $\lambda_{\max}$ and $\lambda_{\min}$ denote the maximum and minimum eigenvalue maps.
Every such matrix $\mtx{A}$ can be expressed as
$$
\mtx{A} = \sum_{i=1}^r \alpha_i \vct{z}_i \vct{z}_i^\adj
\quad\text{where}\quad
\vct{z}_i \in \Z_m^d
\quad\text{and}\quad
m \leq 1 + \frac{1}{2} \sqrt{ (d - 1) \cdot \kappa(\mtx{A}) }.
$$
The coefficients $\alpha_i$ are positive, and the number $r$ of terms
satisfies $r \leq d(d+1)/2$.  The symbol ${}^\adj$ refers to the transpose operation.
\end{theorem}

This result has an alternative interpretation as a matrix factorization:
$$
\mtx{A} = \mtx{Z \Delta Z}^\adj.
$$
In this expression, $\mtx{Z}$ is a $d \times r$ integer matrix with entries bounded by $m$.
The $r \times r$ matrix $\mtx{\Delta}$ is nonnegative and diagonal.

The proof of Theorem~\ref{thm:main} appears in Section~\ref{sec:proof}.
Section~\ref{sec:optimality} demonstrates that the
dependence on the condition number cannot be improved.
We believe that the dependence on the dimension is also
optimal, but we did not find an example that confirms this surmise.

\section{Notation \& Background}

This section contains brief preliminaries.
The books~\cite{HJ90:Matrix-Analysis,Bha97:Matrix-Analysis,Bar02:Course-Convexity,BV04:Convex-Optimization}
are good foundational references for the techniques in this paper.

We use lowercase italic letters, such as $c$, for scalars.
Lowercase boldface letters, such as $\vct{z}$, denote vectors.
Uppercase boldface letters, such as $\mtx{A}$, refer to matrices.
We write $z_i$ for the $i$th component of a vector $\vct{z}$,
and $a_{ij}$ for the $(i, j)$ component of a matrix $\mtx{A}$.
The $j$th column of the matrix $\mtx{A}$ will be denoted by $\vct{a}_j$.

We work primarily in the real linear space $\Sym^d$ of real $d \times d$
symmetric matrices, equipped with the usual componentwise addition and scalar multiplication:
$$
\Sym^d := \big\{ \mtx{A} \in \R^{d \times d} : \mtx{A} = \mtx{A}^\adj \big\}.
$$
Note that $\Sym^d$ has dimension $d(d+1)/2$.
The \term{trace} of a matrix $\mtx{A} \in \Sym^d$ is the sum of its diagonal entries:
$$
\trace(\mtx{A}) := \sum_{i=1}^d a_{ii}.
$$
We equip $\Sym^d$ with the inner product $(\mtx{B}, \mtx{A}) \mapsto \trace(\mtx{BA})$
to obtain a real inner-product space.  All statements about closures refer to the norm
topology induced by this inner product.

Define the set of \term{positive-semidefinite matrices} in $\Sym^d$:
$$
\Sym^d_{+} := \big\{ \mtx{A} \in \Sym^d : \vct{u}^\adj \mtx{A} \vct{u} \geq 0
\text{ for each $\vct{u} \in \R^d$} \big\}.
$$
Similarly, the set of \term{positive-definite matrices} is
$$
\Sym^d_{++} := \big\{ \mtx{A} \in \Sym^d : \vct{u}^\adj \mtx{A} \vct{u} > 0
\text{ for each \emph{nonzero} $\vct{u} \in \R^d$} \big\}.
$$
The members of the set $- \Sym^d_{++}$ are called \term{negative-definite matrices}.

For a matrix $\mtx{A} \in \Sym^d$, the decreasingly ordered eigenvalues will be written as
$$
\lambda^\downarrow_1(\mtx{A}) \geq \lambda^\downarrow_2(\mtx{A})
	\geq \dots \geq \lambda^\downarrow_d(\mtx{A}).
$$
Similarly, the increasingly ordered eigenvalues are denoted as
$$
\lambda^\uparrow_1(\mtx{A}) \leq \lambda^\uparrow_2(\mtx{A})
	\leq \dots \leq \lambda^\uparrow_d(\mtx{A}).
$$
Note that each eigenvalue map $\lambda( \cdot )$ is positively homogeneous; that is,
$\lambda(\alpha \mtx{A}) = \alpha \lambda(\mtx{A})$ for all $\alpha > 0$.

Let us introduce some concepts from conic geometry
in the setting of $\Sym^d$.
A \term{cone} is a subset $K \subset \Sym^d$
that is positively homogeneous; in other words,
$\alpha K = K$ for all $\alpha > 0$.
A \term{convex cone} is a cone that is also a convex set.
The \term{conic hull} of a set $E \subset \Sym^d$
is the smallest convex cone that contains $E$:
\begin{equation} \label{eqn:conic-hull}
\cone(E) := \left\{ \sum_{i=1}^r \alpha_i \mtx{A}_i :
	\alpha_i \geq 0 \text{ and } \mtx{A}_i \in E \text{ and } r \in \N \right\}.
\end{equation}
The conic hull of a finite set is closed.  Since the space $\Sym^d$ has
dimension $d(d+1)/2$, we can choose the explicit value $r = d(d+1)/2$
in the expression~\eqref{eqn:conic-hull}.
This point follows from a careful application of Carath{\'e}odory's
theorem~\cite[Thm.~I(2.3)]{Bar02:Course-Convexity}.

The \term{dual cone} associated with a cone $K \subset \Sym^d$ is the set
\begin{equation} \label{eqn:dual-cone}
K^\adj := \big\{ \mtx{B} \in \Sym^d : \trace(\mtx{BA}) \geq 0 \text{ for each $\mtx{A} \in K$} \big\}.
\end{equation}
This set is always a closed convex cone because it is an intersection of closed halfspaces.
It is easy to check that conic duality reverses inclusion; that is,
for any two cones $C, K \subset \Sym^d$,
$$
C \subset K
\quad\text{implies}\quad
K^\adj \subset C^\adj.
$$
Note that we take the relation $\subset$ to include the possibility that the sets are equal.
The bipolar theorem~\cite[Thm.~IV(4.2)]{Bar02:Course-Convexity} states that the double dual
$(K^\adj)^\adj$ of a cone $K$ equals the closure of the conic hull of $K$.

\section{Proof of Theorem~\ref{thm:main}} \label{sec:proof}

We will establish Theorem~\ref{thm:main} using methods from the geometry
of convex cones.  The result is ultimately a statement about the containment
of one convex cone in another.  We approach this question by verifying
the reverse inclusion for the dual cones.  To obtain a good bound on the
size of the integer vectors, the key idea is to use an averaging argument.

\subsection{Step 1: Reduction to Conic Geometry}

Once and for all, fix the ambient dimension $d$.
First, we introduce the convex cone of positive-definite matrices with bounded condition
number.  For a real number $c \geq 1$, define
$$
K(c) := \big\{ \mtx{A} \in \Sym_{++}^d : \kappa(\mtx{A}) \leq c \big\}.
$$
The set $K(c)$ is a cone because the condition number is scale
invariant: $\kappa(\alpha \mtx{A}) = \kappa(\mtx{A})$ for $\alpha > 0$.
To see that $K(c)$ is convex, write the membership condition
$\kappa(\mtx{A}) \leq c$ in the form
$$
\lambda_{\max}(\mtx{A}) - c \cdot \lambda_{\min}(\mtx{A}) \leq 0.
$$
On the space of symmetric matrices, the maximum eigenvalue is convex, while the minimum eigenvalue is concave~\cite[Ex.~3.10]{BV04:Convex-Optimization}.
Since $K(c)$ is a sublevel set of a convex function, it must be convex.

Next, select a positive integer $m$.
We introduce a closed convex cone of positive-semidefinite matrices
derived from the outer products of bounded integer vectors:
$$
Z(m) := \cone\big\{ \vct{zz}^\adj : \mtx{z} \in \mathbb{Z}^d_m \big\}.
$$
It is evident that $Z(m)$ is a closed convex cone because
it is the conic hull of a finite set.  Note that every element of
this cone can be written as
$$
\sum_{i=1}^r \alpha_i \vct{z}_i \vct{z}_i^\adj
\quad\text{where $\alpha_i \geq 0$ and $\vct{z}_i \in \mathbb{Z}^d_m$.}
$$
By the Carath{\'e}odory Theorem, we may take the number $r$ of summands to be $r = d(d+1)/2$.

Therefore, we can prove Theorem~\ref{thm:main} by verifying that
\begin{equation} \label{eqn:Kc-in-Zm}
K(c) \subset Z(m)
\quad\text{when}\quad
m \geq \frac{1}{2} \sqrt{(d-1) \cdot c}.
\end{equation}
Indeed, the formula $1 + \half \sqrt{(d-1) \cdot \kappa(\mtx{A})}$ in the theorem statement
produces a positive integer that satisfies the latter inequality when $c = \kappa(\mtx{A})$.
Since the operation of conic duality reverses inclusion and $Z(m)$ is closed,
the condition~\eqref{eqn:Kc-in-Zm} is equivalent with
\begin{equation} \label{eqn:Zm*-in-Kc*}
Z(m)^\adj \subset K(c)^\adj
\quad\text{when}\quad
m \geq  \frac{1}{2} \sqrt{(d-1) \cdot c}.
\end{equation}
We will establish the inclusion~\eqref{eqn:Zm*-in-Kc*}.

\subsection{Step 2: The Dual of $K(c)$}

Our next objective is to obtain a formula for the dual cone $K(c)^\adj$.
We claim that
\begin{equation} \label{eqn:Kc-dual}
K(c)^\adj = \left\{ \mtx{B} \in \Sym^d : \sum_{i=s+1}^d \lambda_i^\downarrow(\mtx{B})
	\geq - \frac{1}{c} \sum_{i=1}^s \lambda_i^\downarrow(\mtx{B})
	\text{ where } \lambda^\downarrow_s(\mtx{B}) \geq 0 > \lambda^\downarrow_{s+1}(\mtx{B}) \right\}.
\end{equation}
We instate the convention that $s = d$ when $\mtx{B}$ is positive semidefinite.  In particular,
the set of positive-semidefinite matrices is contained in the dual cone: $\Sym^d_+ \subset K(c)^\adj$.
We also interpret the $s = 0$ case in~\eqref{eqn:Kc-dual} to exclude negative-definite matrices
from $K(c)^\adj$.

Let us establish~\eqref{eqn:Kc-dual}.
The definition~\eqref{eqn:dual-cone} of a dual cone leads to the equivalence
$$
\mtx{B} \in K(c)^\adj
\quad\text{if and only if}\quad
0 \leq \inf_{\mtx{A} \in K(c)} \trace( \mtx{BA} ).
$$
To evaluate the infimum, note that the cone $K(c)$ is orthogonally invariant because the
condition number of a matrix depends only on the eigenvalues.  That is,
$\mtx{A} \in K(c)$ implies that
$\mtx{QAQ}^\adj \in K(c)$ for each orthogonal matrix $\mtx{Q}$ with dimension $d$.
Therefore, $\mtx{B} \in K(c)^\adj$ if and only if
\begin{equation} \label{eqn:Kc-ineq}
0 \leq \inf_{\mtx{A} \in K(c)} \inf_{\mtx{Q}} \trace( \mtx{BQAQ}^\adj )
	= \inf_{\mtx{A} \in K(c)} \sum_{i=1}^d \lambda^\downarrow_i(\mtx{B}) \cdot \lambda^\uparrow_i(\mtx{A}).
\end{equation}
The inner infimum takes place over orthogonal matrices $\mtx{Q}$.
The identity is a well-known result due to Richter~\cite[Satz 1]{Ric58:Abschaetzung};
see the paper~\cite[Thm.~1]{Mir59:Trace-Matrix} for an alternative proof.
This fact is closely related to (a version of)
the Hoffman--Wielandt theorem~\cite[Prob.~III.6.15]{Bha97:Matrix-Analysis}

Now, the members of the cone $K(c)$ are those matrices $\mtx{A}$ whose eigenvalues satisfy
the bounds
$0 < \lambda_{1}^\uparrow(\mtx{A})$ and
$\lambda_{d}^\uparrow(\mtx{A}) \leq c \cdot \lambda^\uparrow_1(\mtx{A})$.
Owing to the invariance of the inequality~\eqref{eqn:Kc-ineq} and the cone $K(c)$ to scaling,
we can normalize $\mtx{A}$ so that $\lambda_1^\uparrow(\mtx{A}) = 1$.
Thus, the inequality~\eqref{eqn:Kc-ineq} holds if and only if
$$
0 \leq \inf\left\{ \sum_{i=1}^d \lambda^\downarrow_i(\mtx{B}) \cdot \mu_i
	: 1 = \mu_1 \leq \mu_2 \leq \dots \leq \mu_d \leq c \right\}.
$$
If $\mtx{B}$ is positive semidefinite, then this bound is always true.
If $\mtx{B}$ is negative definite, then this inequality is always false.
Ruling out these cases, let $s$ be the index where
$\lambda_s^\downarrow(\mtx{B}) \geq 0 > \lambda_{s+1}^\downarrow(\mtx{B})$,
and observe that $0 < s < d$.  The infimum is achieved when we select
$\mu_i = 1$ for $i = 1, \dots s$ and $\mu_i = c$ for $i = s+1, \dots, d$.
In conclusion,
$$
\mtx{B} \in K(c)^\adj
\quad\text{if and only if}\quad
0 \leq \sum_{i=1}^s \lambda_i^\downarrow(\mtx{B}) + c \sum_{i=s+1}^d \lambda_i^\downarrow(\mtx{B}).
$$
With our conventions for $s = 0$ and $s = d$, this inequality coincides with the advertised
result~\eqref{eqn:Kc-dual}.

\subsection{Step 3: The Dual of $Z(m)$}

Next, we check that
\begin{equation} \label{eqn:Zm-dual}
Z(m)^\adj = \big\{ \mtx{B} \in \Sym^d : \vct{z}^\adj \mtx{B} \vct{z} \geq 0
	\text{ for each $\vct{z} \in \Z^d_m$} \big\}.
\end{equation}
According to the definition~\eqref{eqn:dual-cone} of a dual cone,
$$
Z(m)^\adj = \big\{ \mtx{B} \in \Sym^d : \trace(\mtx{BA}) \geq 0
	\text{ for each $\mtx{A} \in Z(m)$} \big\}.
$$
Since $Z(m)$ is the conic hull of the matrices $\vct{zz}^\adj$ where $\vct{z} \in \Z^d_m$,
the matrix $\mtx{B} \in Z(m)^\adj$ if and only if
$\trace(\mtx{BA}) \geq 0$ for each matrix $\mtx{A} = \vct{zz}^\adj$.
Therefore,
$$
Z(m)^\adj = \big\{ \mtx{B} \in \Sym^d : \trace( \mtx{B} \vct{zz}^\adj ) \geq 0
	\text{ for each $\vct{z} \in \Z^d_m$} \big\}.
$$
Cycling the trace, we arrive at the representation~\eqref{eqn:Zm-dual}.

\subsection{Step 4: Checking Membership}

Finally, we need to verify that $Z(m)^\adj \subset K(c)^\adj$ 
under suitable conditions on the parameters $m$ and $c$.

To that end, select a matrix $\mtx{B} \in Z(m)^\adj$.
If $\mtx{B}$ is positive semidefinite, then $\mtx{B} \in K(c)^\adj$ because
$K(c)^\adj$ contains the set of positive-semidefinite matrices.
It is not possible for $\mtx{B}$ to be negative definite because the expression~\eqref{eqn:Zm-dual}
forces $\vct{z}^\adj \mtx{B} \vct{z} \geq 0$ for each nonzero $\vct{z} \in \Z^d_m$.
Therefore, we may exclude these cases.

Let $s$ be the index where $\lambda^\downarrow_s(\mtx{B}) \geq 0 > \lambda^\downarrow_{s+1}(\mtx{B})$,
and note that $0 < s < d$.  The formula~\eqref{eqn:Kc-dual} indicates that we
should examine the sum of the $d - s$ smallest eigenvalues of $\mtx{B}$ to determine
whether $\mtx{B}$ is a member of $K(c)^\adj$.  This sum of eigenvalues can be represented
as a trace~\cite[Eqn.~(4.3.20)]{HJ90:Matrix-Analysis}:
$$
\sum_{i=s+1}^d \lambda^\downarrow_1(\mtx{B})
	= \trace(\mtx{U}^\adj \mtx{B} \mtx{U})
	\quad\text{where $\mtx{U}$ is a $d \times (d-s)$ matrix with orthonormal columns.}
$$
In view of~\eqref{eqn:Zm-dual}, we must use the fact that
$\vct{z}^\adj \mtx{B} \vct{z} \geq 0$ for
$\vct{z} \in \Z^d_m$ to bound the sum of eigenvalues below.

We will achieve this goal with an averaging argument.  For each number $a \in [-1, 1]$,
define an integer-valued random variable:
$$
R_m(a) := \begin{cases}
	\lceil m a \rceil & \text{with probability~$ma - \lfloor ma \rfloor$} \\
	\lfloor m a \rfloor & \text{with probability~$1 - (ma - \lfloor ma \rfloor )$.}
\end{cases}
$$
Each of the random variables $R_m(a)$ is supported on $\{0, \pm 1, \dots, \pm m \}$.
Furthermore, $\Expect R_m(a) = ma$ and $\Var(R_m(a)) \leq \frac{1}{4}$.
In other words, we randomly round $ma$ up or down to the nearest integer
in such a way that the average value is $ma$ and the variance is uniformly bounded.
Note that $R_m(a)$ is a constant random variable whenever $ma$ takes an integer value.

We apply this randomized rounding operation to each entry $u_{ij}$ of the matrix $\mtx{U}$.
Let $\mtx{X}$ be a $d \times (d-s)$ random matrix with independent entries $X_{ij}$
that have the distributions
$$
X_{ij} \sim \frac{1}{m} R_m(u_{ij})
\quad\text{for $i = 1, \dots, d$ and $j = 1, \dots, d-s$.}
$$
By construction, $\Expect \mtx{X} = \mtx{U}$ and
$\Var(X_{ij}) \leq 1/(4m^2)$ for each pair $(i, j)$ of indices.

Develop the quantity of interest by adding and subtracting the random matrix $\mtx{X}$:
$$
\sum_{i=s+1}^d \lambda^\downarrow_1(\mtx{B})
	= \trace(\mtx{U}^\adj \mtx{B} \mtx{U})
	= \trace(\mtx{X}^\adj \mtx{B} \mtx{X})
	- \trace((\mtx{X} - \mtx{U})^\adj \mtx{B}(\mtx{X} - \mtx{U}))
	- \trace(\mtx{U}^\adj \mtx{B} (\mtx{X} - \mtx{U}))
	- \trace((\mtx{X} - \mtx{U})^\adj \mtx{B} \mtx{U}).
$$
Take the expectation over $\mtx{X}$ and use the property $\Expect \mtx{X} = \mtx{U}$
to reach
\begin{equation} \label{eqn:eig-avg}
\sum_{i=s+1}^d \lambda^\downarrow_1(\mtx{B})
	= \Expect \trace(\mtx{X}^\adj \mtx{B} \mtx{X})
	- \Expect \trace((\mtx{X} - \mtx{U})^\adj \mtx{B}(\mtx{X} - \mtx{U})).
\end{equation}
It remains to bound the right-hand side of~\eqref{eqn:eig-avg} below.

Expand the trace in the first term on the right-hand side of~\eqref{eqn:eig-avg}:
\begin{equation} \label{eqn:sum-eigs-1}
\Expect \trace( \mtx{X}^\adj \mtx{B} \mtx{X} )
	= \Expect \left[ \sum_{j=1}^{d-s} \vct{x}_{j}^\adj \mtx{B} \vct{x}_{j} \right]
	= \frac{1}{m^2} \Expect \left[ \sum_{j=1}^{d-s} (m\vct{x}_{j})^\adj \mtx{B} (m\vct{x}_{j}) \right]
	\geq 0.
\end{equation}
We have written $\vct{x}_{j}$ for the $j$th column of $\mtx{X}$.
Each vector $m \vct{x}_{j}$ belongs to $\Z^d_m$.  Since
$\mtx{B} \in Z(m)^\adj$, it follows from the representation~\eqref{eqn:Zm-dual}
of the cone that each of the summands is nonnegative.

Next, we turn to the second term on the right-hand side of~\eqref{eqn:eig-avg}.
$$
\Expect \trace((\mtx{X} - \mtx{U})^\adj \mtx{B} (\mtx{X} - \mtx{U}))
	= \sum_{j=1}^{d-s} \Expect\big[ (\vct{x}_{j} - \vct{u}_{j})^\adj \mtx{B} (\vct{x}_{j} - \vct{u}_{j}) \big]
	= \sum_{j=1}^{d-s} \sum_{i=1}^d \Expect\big[ (X_{ij} - u_{ij})^2 \big] \cdot b_{ii}
	\leq \frac{d-s}{4m^2} \sum_{i=1}^d (b_{ii})_+.
$$
In the second identity, we applied the fact that the entries of the vector $\vct{x}_{j} - \vct{u}_{j}$
are independent, centered random variables to see that there is no contribution from the off-diagonal
terms of $\mtx{B}$.  The inequality relies on the variance bound $1/(4m^2)$ for each random variable $X_{ij}$.  
The function $(\cdot)_+ : a \mapsto \max\{a, 0\}$ returns the positive part of a number.

Schur's theorem~\cite[Thm.~4.3.26]{HJ90:Matrix-Analysis} states that eigenvalues of the
symmetric matrix $\mtx{B}$ majorize its diagonal entries.
Since $(\cdot)_+$ is convex, the real-valued map $\vct{a} \mapsto \sum_{i=1}^d (a_i)_+$
on $\R^d$ respects the majorization relation~\cite[Thm.~II.3.1]{Bha97:Matrix-Analysis}.  Thus,
$$
\sum_{i=1}^d (b_{ii})_+ \leq \sum_{i=1}^d (\lambda^\downarrow_i(\mtx{B}))_+
	= \sum_{i=1}^s \lambda^\downarrow_i(\mtx{B}).
$$
The equality relies on the assumption that the eigenvalues $\lambda_i^\downarrow(\mtx{B})$ become negative
at index $s+1$.

Merging the last two displays, we obtain the estimate
\begin{equation} \label{eqn:sum-eigs-2}
\Expect \trace((\mtx{X} - \mtx{U})^\adj \mtx{B} (\mtx{X} - \mtx{U}))
	\leq \frac{d-s}{4m^2} \sum_{i=1}^s \lambda^\downarrow_i(\mtx{B}).
\end{equation}
This bound has exactly the form that we need.

Combining~\eqref{eqn:eig-avg}, \eqref{eqn:sum-eigs-1}, and \eqref{eqn:sum-eigs-2},
we arrive at the inequality
$$
\sum_{i=s+1}^d \lambda^\downarrow_1(\mtx{B})
	\geq - \frac{d-s}{4m^2} \sum_{i=1}^s \lambda^\downarrow_i(\mtx{B}).
$$
In view of the representation~\eqref{eqn:Kc-dual} of the dual cone $K(c)^\adj$,
the matrix $\mtx{B} \in K(c)^\adj$ provided that
$$
- \frac{d-s}{4m^2} \geq - \frac{1}{c}.
$$
Rearranging this expression, we obtain the sufficient condition
$$
m \geq \frac{1}{2} \sqrt{(d - s) \cdot c}
\quad\text{implies}\quad
\mtx{B} \in K(c)^\adj.
$$
For a general matrix $\mtx{B} \in Z(m)^\adj$,
we do not control the index $s$ where the eigenvalues of $\mtx{B}$ change sign,
so we must insulate ourselves against the worst case, $s = 1$.
This choice leads to the condition~\eqref{eqn:Zm*-in-Kc*},
and the proof is complete.

\section{Optimality} \label{sec:optimality}

There are specific matrices where the size of the integers
in the representation does not depend on the condition number.
For instance, let $b \geq 1$, and consider the matrix
$$
\mtx{A} = \begin{bmatrix} b & 0 \\ 0 & 1 \end{bmatrix}
	= b \begin{bmatrix} 1 \\ 0 \end{bmatrix} \begin{bmatrix} 1 \\ 0 \end{bmatrix}^\adj
	+ 1 \begin{bmatrix} 0 \\ 1 \end{bmatrix} \begin{bmatrix} 0 \\ 1 \end{bmatrix}^\adj.
$$
The condition number $\kappa(\mtx{A}) = b$, which we can make arbitrarily large,
but the integers in the representation never exceed one.

Nevertheless, we can show by example that the dependence of Theorem~\ref{thm:main}
on the condition number is optimal in dimension $d = 2$.
For a number $b \geq 1$, consider the $2 \times 2$ matrix
$$
\mtx{A} = \begin{bmatrix} b^2 + 1 & b \\ b & 2 \end{bmatrix}
	= \begin{bmatrix} b \\ 1 \end{bmatrix} \begin{bmatrix} b \\ 1 \end{bmatrix}^\adj
	+ \begin{bmatrix} 1 & 0 \\ 0 & 1 \end{bmatrix}.
$$
From this representation, we quickly determine that the eigenvalues of $\mtx{A}$
are $1$ and $b^2 + 2$, so the condition number $\kappa(\mtx{A}) = b^2 + 2$.

Suppose that we can represent the matrix $\mtx{A}$ as a positive linear combination
of outer products of vectors in $\Z_m^2$.  We need at most $d(d+1)/2 = 3$ summands:
\begin{equation} \label{eqn:A-example-rep}
\mtx{A} = \alpha \begin{bmatrix} x_1 \\ x_2 \end{bmatrix} \begin{bmatrix} x_1 \\ x_2 \end{bmatrix}^\adj
	+ \beta \begin{bmatrix} y_1 \\ y_2 \end{bmatrix} \begin{bmatrix} y_1 \\ y_2 \end{bmatrix}^\adj
	+ \gamma \begin{bmatrix} z_1 \\ z_2 \end{bmatrix} \begin{bmatrix} z_1 \\ z_2 \end{bmatrix}^\adj
	\quad\text{where $\alpha, \beta, \gamma > 0$ and $x_i, y_i, z_i \in \Z_m^1$.}
\end{equation}
The equations in~\eqref{eqn:A-example-rep}
associated with the top-left and bottom-right entries of $\mtx{A}$ read as
\begin{equation} \label{eqn:A-example-eqns}
b^2 + 1 = \alpha x_1^2 + \beta y_1^2 + \gamma z_1^2
\quad\text{and}\quad
2 = \alpha x_2^2 + \beta y_2^2 + \gamma z_2^2.
\end{equation}
We consider thee cases: (i) all three of $x_2, y_2, z_2$ are nonzero;
(ii) exactly two of $x_2, y_2, z_2$ are nonzero; and
(iii) exactly one of $x_2, y_2, z_2$ is nonzero.

Let us begin with case (i).  Since $x_2$, $y_2$, and $z_2$ take nonzero integer values,
the second equation in~\eqref{eqn:A-example-eqns} ensures that
$$
2 \geq (\alpha + \beta + \gamma) \min\big\{x_2^2, y_2^2, z_2^2 \big\} \geq \alpha + \beta + \gamma.
$$
Introducing this fact into the first equation in~\eqref{eqn:A-example-eqns}, we find that
$$
b^2 + 1	\leq (\alpha + \beta + \gamma) \max\big\{ x_1^2, y_1^2, z_1^2 \big\}
	\leq 2 \max\big\{ x_1^2, y_1^2, z_1^2 \big\}.
$$
We obtain a lower bound on the magnitude $m$ of integers in a representation
of $\mtx{A}$ where $x_2, y_2, z_2$ are all nonzero:
\begin{equation} \label{eqn:case-i}
m \geq \max\big\{ \abs{x_1}, \abs{y_1}, \abs{z_1} \big\}
	\geq \frac{1}{\sqrt{2}} \sqrt{b^2 + 1}
	= \frac{1}{\sqrt{2}} \sqrt{\kappa(\mtx{A}) - 1}.
\end{equation}
Since the bound~\eqref{eqn:case-i} is worse than the estimate in
Theorem~\ref{thm:main} for large $b$, we discover that the optimal
integer representation of $\mtx{A}$ has at least one zero among $x_2, y_2, z_2$.

Next, we turn to case (ii).  By symmetry, we may assume that $z_2 = 0$.
As before, the second equation in~\eqref{eqn:A-example-eqns} shows that
$\alpha + \beta \leq 2$.  Meanwhile,
the representation~\eqref{eqn:A-example-rep} implies that
$$
\mtx{A} - \gamma \begin{bmatrix} z_1^2 & 0 \\ 0 & 0 \end{bmatrix}
	= \begin{bmatrix} b^2 + 1 - \gamma z_1^2 & b \\ b & 2 \end{bmatrix}
	\quad\text{is positive semidefinite.}
$$
Since the determinant of a positive-semidefinite matrix is nonnegative,
we find that $0 \leq 2(b^2 + 1 - \gamma z_1^2) - b^2$.
Equivalently, $\gamma z_1^2 \leq \frac{1}{2}(b^2 + 2)$.  The first equation in~\eqref{eqn:A-example-eqns}
now delivers
$$
b^2 + 1 = \alpha x_1^2 + \beta y_1^2 + \gamma z_1^2
	\leq 2 \max\big\{ x_1^2, y_1^2 \big\} + \frac{1}{2} \big( b^2 + 2 \big).
$$
It follows that $\max\{ x_1^2, y_1^2 \} \geq b^2 / 4$.  We obtain a lower bound
on the magnitude $m$ of the integers in a representation of $\mtx{A}$ where
two of $x_2, y_2, z_2$ are nonzero:
\begin{equation} \label{eqn:case-ii}
m \geq \max\big\{ \abs{x_1}, \abs{y_1} \big\} \geq \frac{1}{2} b = \frac{1}{2} \sqrt{\kappa(\mtx{A}) - 2}.
\end{equation}
In case (iii), a similar argument leads to the same lower bound for $m$.

Examining~\eqref{eqn:case-ii}, we surmise that the bound from Theorem~\ref{thm:main}
$$
m \leq 1 + \frac{1}{2} \sqrt{(d - 1) \cdot \kappa(\mtx{A})},
$$
on the magnitude $m$ of integers in a representation of $\mtx{A}$
cannot be improved when $d = 2$ and the condition number $\kappa(\mtx{A})$
becomes large.  Considering the $d\times d$ matrix
$$
\begin{bmatrix} \mtx{A} & \mtx{0} \\ \mtx{0} & \Id_{d-2}
\end{bmatrix},
$$
an analogous arguments proves that the dependence of Theorem~\ref{thm:main}
on the condition number is optimal in every dimension $d$.

\section*{Acknowledgments}

The author wishes to thank Andrew Barbour for calling this problem to his attention
and to Madeleine Udell for describing the application to analog-to-digital matrix multiplication.
This research was undertaken at the Institute for Mathematical Sciences (IMS)
at the National University of Singapore (NUS) during the workshop
on New Directions in Stein's Method in May 2015.
The author gratefully acknowledges
support from ONR award N00014-11-1002 and the Gordon \& Betty Moore Foundation.

\bibliographystyle{myalpha}

\begin{thebibliography}{LUW15}

\bibitem[Bar02]{Bar02:Course-Convexity}
A.~Barvinok.
\newblock {\em A course in convexity}, volume~54 of {\em Graduate Studies in
  Mathematics}.
\newblock American Mathematical Society, Providence, RI, 2002.

\bibitem[Bha97]{Bha97:Matrix-Analysis}
R.~Bhatia.
\newblock {\em Matrix analysis}, volume 169 of {\em Graduate Texts in
  Mathematics}.
\newblock Springer-Verlag, New York, 1997.

\bibitem[BLX15]{BLX15:Multivariate-Approximation}
A.~D. Barbour, M.~J. Luczak, and A.~Xia.
\newblock Multivariate approximation in total variation.
\newblock Manuscript, May 2015.

\bibitem[BV04]{BV04:Convex-Optimization}
S.~Boyd and L.~Vandenberghe.
\newblock {\em Convex optimization}.
\newblock Cambridge University Press, Cambridge, 2004.

\bibitem[HJ90]{HJ90:Matrix-Analysis}
R.~A. Horn and C.~R. Johnson.
\newblock {\em Matrix analysis}.
\newblock Cambridge University Press, Cambridge, 1990.
\newblock Corrected reprint of the 1985 original.

\bibitem[LUW15]{LUW15:Factorization-Analog}
E.~H. Lee, M.~Udell, and S.~S. Wong.
\newblock Factorization for analog-to-digital matrix multiplication.
\newblock In {\em Proc. 40th Intl. Conf. Acoustics, Speech, and Signal
  Processing (ICASSP)}, Brisbane, 2015.

\bibitem[Mir59]{Mir59:Trace-Matrix}
L.~Mirsky.
\newblock On the trace of matrix products.
\newblock {\em Math. Nachr.}, 20:171--174, 1959.

\bibitem[Ric58]{Ric58:Abschaetzung}
H.~Richter.
\newblock Zur {A}bsch\"atzung von {M}atrizennormen.
\newblock {\em Math. Nachr.}, 18:178--187, 1958.

\end{thebibliography}

\end{document}